\documentclass[10pt]{amsart}
\usepackage[centertags]{amsmath}
\usepackage{graphicx}
\usepackage{amsfonts,amssymb,latexsym,epsfig,amsthm}
\usepackage{fullpage}

\usepackage{color}
\definecolor{blue}{rgb}{0,0,1}
\definecolor{red}{rgb}{1,0,0}

\newtheorem{theorem}{Theorem}
\newtheorem{lemma}{Lemma}[section]
\newtheorem{corollary}[lemma]{Corollary}

\newcommand{\EE}{\mathbb{E}}
\newcommand{\NN}{\mathbb{N}}
\newcommand{\PP}{\mathbb{P}}
\newcommand{\RR}{\mathbb{R}}
\newcommand{\BT}{\tilde{B}}

\newcommand{\Go}{\overline{G}_k}
\newcommand{\Gu}{\underline{G}_k}

\title{A GUE Central Limit Theorem and Universality of Directed First and Last Passage Site Percolation}

\begin{document}
\author{Jinho Baik  \and Toufic M. Suidan}
\address{Jinho Baik (baik@umich.edu): Department of Mathematics,
  University of Michigan,  Ann Arbor, MI}
\address{Toufic Suidan (suidan@cims.nyu.edu): Department of Mathematics, University of
  California, Santa Cruz, CA, and Courant
  Institute of Mathematical Sciences, New York University, NYC, NY}

\begin{abstract}
  We prove a GUE central limit theorem for random variables
  with finite fourth moment. We apply this theorem to prove that the
  directed first and last passage percolation problems in thin
  rectangles exhibit universal fluctuations given by the Tracy-Widom
  law.
\end{abstract}
\maketitle

\section{Introduction}\label{sec:intro}

In the last few years, it has become clear that random matrix
theory is intimately related to a variety of questions arising in
physics, statistics, combinatorics, representation theory, number
theory, and probability theory (see e.g. \cite{Mehta, Deiftbook,
KatzSarnak, TWICM}). It is the fluctuations of random matrix
ensembles which are common to the specific problems from each of
these areas. In the hope of understanding the universal nature of
the random matrix distributions, it is natural to search for
central limit theorems for which these distributions are the
limiting objects.

The following is a Gaussian Unitary Ensemble (GUE) central limit
theorem.

\begin{theorem}[GUECLT]\label{maintheorem}
Suppose that $\{X_i^j\}_{i,j=1}^\infty$ is a family of independent identically
distributed random variables such that $\EE X_i^j=0$, $\EE |X_i^j|^2=1$,
and $\EE |X_i^j|^4<\infty$. Let
\begin{eqnarray}
L(N,k)=\sup_{0=i_0\leq i_1\leq...\leq i_k=N}
\sum_{j=1}^{k}\bigl\{\sum_{i=i_{j-1}+1}^{i_j} X_i^j \bigr\}, \\
R(N,k)=\inf_{0=i_0\leq i_1\leq...\leq i_k=N}
\sum_{j=1}^{k}\bigl\{\sum_{i=i_{j-1}+1}^{i_j} X_i^j \bigr\}.
\end{eqnarray}
If $k,N\rightarrow \infty$ such that $k=o(N^\alpha)$, $\alpha<\frac{3}{14}$, then
\begin{eqnarray}\label{limit}
\Bigl(\frac{L(N,k)}{N^\frac{1}{2}}-2\sqrt{k}\Bigr) k^\frac{1}{6}
\Longrightarrow F_{GUE},\\
\label{limit2} \Bigl(-\frac{R(N,k)}{N^\frac{1}{2}}-2\sqrt{k}\Bigr)
k^\frac{1}{6} \Longrightarrow F_{GUE},
\end{eqnarray}
where $F_{GUE}$ is the GUE Tracy-Widom distribution. If
$\{X_i^j\}_{i,j=1}^\infty$ are independent identically distributed
Gaussian random variables satisfying the above conditions,
then~\eqref{limit} and~\eqref{limit2} hold if
$k=o(N^\alpha)$, $\alpha<\frac{3}{7}$.
\end{theorem}

The GUE Tracy-Widom distribution function \cite{TW1} is given by
\begin{equation}
  F_{GUE}(x)= \exp\biggl\{ -\int_x^\infty (s-x)q^2(s)ds \biggr\},
\end{equation}
where $q(x)$ solves the Painlev\'e II equation,
\begin{equation}
  q''=2q^3+xq
\end{equation}
subject to the condition that $q(x) \sim Ai(x)$ as $x\to +\infty$;
$Ai(x)$ denotes the Airy function. The function $F_{GUE}$ is
the limiting distribution function for the largest eigenvalue of
the Gaussian Unitary Ensemble as the dimension of the matrices grows to
infinity (see~\eqref{eq:GUElargest}).

Theorem~\ref{maintheorem} is intimately related to the directed
first and last passage percolation problems. Consider the $\NN \times \NN$ lattice and a set of associated
independent identically distributed random variables
$\{X_i^j\}_{i,j=1}^\infty$ satisfying $\EE X_i^j =\mu$,
$\EE|X_i^j|^2-\mu^2=\sigma^2$, and $\EE |X_i^j|^4 <\infty$. An
\emph{up/right path} $\pi$ from the site $(1,1)$ to the site
$(N,k)$ is a collection of sites $\{(i_k, j_k)\}_{k=1}^{N+k-1}$
satisfying $(i_1, j_1)=(1,1)$, $(i_{N+k-1}, j_{N+k-1})=(N,k)$ and
$(i_{k+1}, j_{k+1})-(i_k,j_k)$ is either $(1,0)$ or
$(0,1)$. Let $(1,1)\nearrow (N,k)$ denote the set of such up/right
paths. The directed first and last passage times to $(N,k)\in \NN
\times \NN$, denoted by $L^f(N,k)$ and $L^l(N,k)$, respectively,
are defined as
\begin{eqnarray}
   L^f(N,k)&=& \min_{\pi\in (1,1)\nearrow (N,k)} \sum_{(i,j)\in
   \pi} X_i^j ,\\
   L^l(N,k)&=& \max_{\pi\in (1,1)\nearrow (N,k)} \sum_{(i,j)\in
   \pi} X_i^j.
\end{eqnarray}
If $X_i^j$ is interpreted as the time to pass the site $(i,j)$,
$L^f(N,k)$ and $L^l(N,k)$ represent the minimal and the maximal
time to travel from the site $(1,1)$ to $(N,k)$ along an
admissible path. Since the directed last passage percolation time
can be viewed as the departure time in queuing theory (see e.g.
\cite{GlynnWhitt}), the result below also applies to
queuing theory. In addition, Corollary~\ref{LPP} also applies to the
flux of particles at a given site in the totally asymmetric simple
exclusion process (see e.g.~\cite{Seppal1997}). Theorem~\ref{maintheorem} implies the following.

\begin{corollary}\label{LPP}
Suppose that $\{X_i^j\}_{i,j=1}^\infty$ is a family of independent
identically distributed random variables such that $\EE
X_i^j=\mu$, $\EE |X_i^j|^2-\mu^2=\sigma^2$, and $\EE |X_i^j|^4<\infty$. For any $s\in \mathbb{R}$,
\begin{eqnarray}\label{eq:l}
\lim_{N,k\rightarrow\infty} \PP\biggl( \frac{L^l(N,k)-\mu
(N+k-1)- 2\sigma\sqrt{Nk}}{\sigma k^{-1/6}N^\frac{1}{2}} \leq s \biggr) &=& F_{GUE}(s), \\
\label{eq:f} \lim_{N,k\rightarrow\infty}\PP\biggl(
\frac{L^f(N,k)-\mu (N+k-1)+2\sigma\sqrt{Nk}}{\sigma
k^{-1/6}N^\frac{1}{2}}  \leq s \biggr) &=& 1-F_{GUE}(-s),
\end{eqnarray}
where $k=o(N^\alpha)$ and $\alpha<\frac{3}{14}$. If $\{X_i^j\}_{i,j=1}^\infty$
are independent identically distributed Gaussian random variables
satisfying the above conditions, then~\eqref{eq:l}
and~\eqref{eq:f} hold when $k=o(N^\alpha)$, $\alpha<\frac{3}{7}$.
\end{corollary}

The restriction $k=o(N^\alpha)$, $\alpha< \frac{3}{14}$, in the above corollary seems to be a technical
matter. It is believed that such a central limit theorem (after
possible changes in scaling) should hold as $N,k \to\infty$ with no
restrictions on their relative rate of growth. When $N$ and $k$
are of the same order, such a limit law for $L^l$ was obtained for
geometric and exponential random variables by Johansson
\cite{kurtj:shape}\footnote{For slightly different choices of admissible
paths, a few more explicit random variables are shown to have the same type of limit law (see e.g. \cite{BDJ,
GTW})} using Young tableaux theory in combinatorics and techniques
from random matrix theory. However, analyzing the last passage
percolation time for general random variables in arbitrary scaling
regimes is still an open question.

Even the determination of the `time constants' as $N, k\to\infty$
for general random variables is a challenge. Glynn and Whitt \cite{GlynnWhitt} proved that for $k=o(N)$,
\begin{equation}\label{eq:GlynnWhitt}
  \lim_{N,k\to\infty} \frac{L^l(N,k)-\mu N}{\sqrt{Nk}} = \alpha,
\end{equation}
$\PP$-almost surely for some $\alpha$, where $\alpha$ was proved
to be equal to $2\sigma$ by Sepp\"al\"ainen \cite{Seppal1997} (see
also~\cite{Martin}). Note that the definition of the last passage
percolation, $L^l(N,k)$, is symmetric in $N$ and $k$.
Hence,~\eqref{eq:GlynnWhitt} can not be true when $N$ and $k$ are
of the same order. For exponential random variables
of mean $1$, an interpretation of a theorem of Rost
\cite{Rost1981} yields that
\begin{equation}\label{eq:exp}
  \lim_{n\to\infty} \frac{L^l([xn], [yn])}{n} = (\sqrt{x}+\sqrt{y})^2.
\end{equation}
Johansson~\cite{kurtj:shape} (see equation (1.4)) proved that for geometric random variables of parameter $q$, 
\begin{equation}\label{eq:geom}
  \lim_{n\to\infty} \frac{L^l([xn], [yn])}{n} =
\frac{q(x+y)+2\sqrt{qxy}}{1-q}.
\end{equation}
These two cases seem to be the only explicitly computed examples when
$N$ and $k$ are of the same order.

\bigskip

Note that the directed first and last passage times can also be
written as
\begin{eqnarray}
L^f(N,k)&=& \inf_{1=i_0\leq i_1\leq...\leq i_k=N}
\sum_{j=1}^{k}\biggl\{\sum_{i=i_{j-1}}^{i_j} X_i^j \biggr\},\\
L^l(N,k)&=& \sup_{1=i_0\leq i_1\leq...\leq i_k=N}
\sum_{j=1}^{k}\biggl\{\sum_{i=i_{j-1}}^{i_j} X_i^j \biggr\}.
\end{eqnarray}
When $k$ is \emph{fixed} and
$N\to\infty$, the last passage time is obtained along a path that
lies on the first level, $j=1$, for a certain `time', $\lfloor t_1N\rfloor$, then jumps to the
second level, $j=2$, and stays there until the `time' $\lfloor
t_2N\rfloor$, and so on such that the total sum is maximized. As the sum of random variables on each
level is `asymptotically equal' to a Brownian motion after proper
centering and scaling, by the Donsker's theorem one can imagine that
\begin{equation}
  \frac{L^l(N,k)-\mu N}{\sigma\sqrt{N}} \Longrightarrow
  \hat{D}_k(1) := \sup_{0=t_0\le t_1\le \dots\le t_k=1}
  \sum_{j=1}^k \bigl( B_j(t_j)-B_j(t_{j-1}) \bigr),
\end{equation}
where $B_j(t)$, $j=1,\dots, k$, are independent Brownian motions.
This was proved in \cite{GlynnWhitt} in the context of queuing
theory.

The distribution function of $\hat{D}_k(1)$ is difficult to compute since
$\hat{D}_k(1)$ is a complicated and unusual functional of $k$
Brownian motions. However, as mentioned earlier, for some special random
variables, Young tableaux theory in combinatorics can be applied
to explicitly compute the distribution of $L^l$. By
studying a special random variable, Baryshnikov
\cite{Baryshnikov} and Gravner, Tracy and Widom \cite{GTW} (see
also \cite{OConnellYor2002, OConnell2003}) showed that
\begin{equation}
  \PP(\hat{D}_k(1) \le s) = \frac{1}{Z_k} \int_{-\infty}^s \dots
  \int_{-\infty}^s \prod_{1\le i\le j\le k} |\xi_i-\xi_j|^2
  \prod_{j=1}^k e^{-\frac12 \xi_j^2} d\xi_j,
\end{equation}
where $Z_k=1!2!\dots k!(2\pi)^{k/2}$. As the integrand on the
right-hand-side is the density of the eigenvalues $\xi_1, \dots,
\xi_k$ of $k\times k$ GUE random matrix, this identity implies
that $\hat{D}_k(1)$ has precisely the same distribution as the
largest eigenvalue, $\xi_{\max}(k)$, of $k\times k$ GUE random
matrix. It is well-known in random matrix theory that (see e.g.
\cite{Forrester, TW1})
\begin{equation}\label{eq:GUElargest}
  \lim_{k\to\infty} \PP \bigl( \bigl(\xi_{\max}(k) -2\sqrt{k }\bigr) k^{1/6} \le s \bigr)
  = F_{GUE}(s).
\end{equation}
Therefore,
\begin{equation}
  \lim_{k\to\infty} \lim_{N\to\infty} \PP \biggl(
  \biggl\{\frac{L^l(N,k)-\mu N}{\sigma\sqrt{N}} - 2\sqrt{k}
  \biggr\}k^{1/6} \le s \biggr) = F_{GUE}(s).
\end{equation}
The content of Corollary~\ref{LPP} is that the two limits can be
taken simultaneously as long as $k=o(N^{\alpha})$, $\alpha< \frac{3}{14}$.
Note that in the centering $\mu(N+k-1)$ of~\eqref{eq:l}, the term
$\mu (k-1)$ accounts for the `up movements' of the path while $\mu
N$ comes from the `right movements' of the path. The term $\mu
(k-1)$ is negligible when $k=o(N^{\frac{3}{7}})$, but it is retained
in the formula in order to emphasize the symmetry of $N,k$. \bigskip

Theorem~\ref{maintheorem} is a central limit theorem whose limit is
the Tracy-Widom GUE top and bottom eigenvalue distributions. The joint
distribution of the top and bottom $n$ eigenvalues is also a natural limit of a
central limit theorem which depends on the O'Connell and Yor
representation of the $k$-dimensional Dyson
process~\cite{OConnellYor2002,
  OConnell2003}. Following~\cite{OConnellYor2002, OConnell2003}, let
$C_0([0,1], \mathbb{R})$ be the space of continuous functions
$f:[0,1]\rightarrow \mathbb{R}$ with $f(0)=0$. For $f,g\in C_0([0,1],\mathbb{R})$, define
\begin{eqnarray}
f\otimes g(t)=\inf_{0\leq s\leq t} [f(s)+g(t)-g(s)], \\
f\odot g(t)=\sup_{0\leq s\leq t} [f(s)+g(t)-g(s)].
\end{eqnarray}
The order of operations is from left to right. Define a sequence of
mappings $\Gamma_k:C_0([0,1],\mathbb{R})^k \rightarrow C_0([0,1],\mathbb{R})^k$ by
\begin{equation}
\Gamma_2(f,g)=(f\otimes g, g\odot f),
\end{equation}
and for $k>2$,
\begin{eqnarray}
\Gamma_k(f_1,...,f_k)&=&(f_1\otimes f_2 \otimes f_3\otimes
\cdots \otimes f_k, \\ &&\Gamma_{k-1}(f_2\odot f_1, f_3\odot
(f_1\otimes f_2), ..., f_k\odot (f_1 \otimes \cdots
\otimes f_{k-1}))).\nonumber
\end{eqnarray}
O'Connell and Yor proved that if $B_1,...,B_k$ are independent
standard one-dimensional Brownian motions, $\Gamma_k(B_1,...,B_k)$ has the same distribution on path space,
$C_0([0,1],\mathbb{R})^k$, as $k$ one-dimensional Brownian motions starting from
the origin, conditioned (in the sense of Doob) never to collide. It is
well known that this process can be interpreted as the $k$-dimensional
GUE Dyson eigenvalue process. The first coordinate of $\Gamma_k$ is
the smallest eigenvalue, the second coordinate of $\Gamma_k$ is the second smallest
eigenvalue, and so on. The proof of Theorem~\ref{maintheorem}
immediately implies the following theorem once one notes that the
first $n$ coordinates of $\Gamma_k$ are
Lipshitz functions with a fixed Lipshitz constant independent of
$k$. We omit the details.
\begin{theorem} Suppose that $\{X_i^j\}_{i,j=1}^\infty$ is a family of
  independent identically distributed random variables such that $\EE
  X_i^j=0$, $\EE |X_i^j|^2=1$, and $\EE |X_i^j|^4<\infty$. Let
  $s_N^j(i)=\frac{1}{\sqrt{N}} \sum_{l=1}^i X_l^j$. As $N,k\rightarrow
  \infty$ such that $k=O(N^\alpha)$, $\alpha<\frac{3}{14}$,
\begin{eqnarray}
(\Gamma_k(s_N^1,...,s_N^k)_1(N),...,\Gamma_k(s_N^1,...,s_N^k)_n(N))
\Longrightarrow F_{GUE}^n,
\end{eqnarray}
where $F_{GUE}^n$ denotes the limiting joint
probability distribution of the bottom $n$ eigenvalues of the GUE. By
symmetry, a similar statement holds for the top $n$ eigenvalues.
\end{theorem}



Proving a GUE central limit theorem and analyzing certain scaling regimes
for the directed first and last passage percolation problems is the
main focus of this paper. It is interesting to seek central limit type
theorems for the Gaussian Orthogonal Ensemble (GOE) and the Gaussian
Symplectic Ensemble (GSE) Tracy-Widom distributions. The arguments of
this paper can be adapted to prove such theorems. The authors will
address this topic in several forthcoming papers. The authors will
also describe the fluctuations of restricted Plancharel measures from
this point of view.

\medskip
\noindent{\bf Note.}
While completing this paper, the authors learned that T. Bodineau and
J. Martin had announced a similar
result~\cite{BodineauMartin}. Both~\cite{BodineauMartin} and this
paper use the idea of approximating random walks by Brownian motion in
order to analyze the last passage time. As opposed
to~\cite{BodineauMartin}, the authors use the Skorohod embedding
theorem to couple random walks to Brownian motion. Bodineau and Martin
use the Koml\'os, Major and Tusn\'ady approximation theorem (which
produces tighter bounds) to achieve such a coupling. The authors are
grateful to Rongfeng Sun for bringing~\cite{BodineauMartin} to their
attention.

\medskip
\noindent {\bf Acknowledgments.} The authors thank S.R.S. Varadhan for
his crucial comments at an early stage of this work. The authors thank
G. Ben Arous, P. Deift, and R. Sun for useful discussions. The authors
also thank J. Martin for constructive comments on the first
version of this paper. The work of Baik was supported in part by NSF
Grant \#DMS-0350729 and the AMS Centennial Fellowship. The work of
Suidan was supported in part by a NSF postdoctoral fellowship.

\section{Proofs}\label{sec:proof}

Theorem~\ref{maintheorem} and Corollary~\ref{LPP} are proven in this section.

\begin{proof}[Proof of Theorem~\ref{maintheorem}]

The proof is based on embedding general random walks in Brownian
motions by means of the classical Skorohod embedding theorem
\cite{Billingsley}.
\begin{theorem}[Skorohod Embedding Theorem] \label{Skorohod}
Let $\{B_t\}_{t\in \RR^+}$ be a one-dimensional standard Brownian
motion and $X$ be a real valued random variable satisfying $\EE X=0$, $\EE
X^2=1$, and $\EE X^4 < \infty$. Then, there exists a Brownian stopping
time, $\tau$, such that $B_\tau$ is distributed as $X$, $\EE \tau=\EE
X^2=1$, and $\EE \tau^2 \leq 4\EE X^4$.
\end{theorem}
Let $\{X_i\}_{i=1}^\infty$ be a sequence of independent identically
distributed random variables satisfying the assumptions of
Theorem~\ref{Skorohod}. There exists a sequence of independent
identically distributed positive random variables $\{\tau_i
\}_{i=1}^\infty$ such that $B_{\tau_1+...+\tau_n}$ is distributed as
$X_1+...+X_n$ and
$\{B_{\tau_1+...+\tau_{n+1}}-B_{\tau_1+...+\tau_n}\}_{n=0}^\infty$ is
a sequence of independent identically distributed random variables
with distribution $X_1$.  This follows from Theorem~\ref{Skorohod} and
the strong Markov property of Brownian motion.

To prove uniform estimates in time for the difference of embedded
random walks and the scaled Brownian motion in which the walks are
embedded, the following two real valued processes are useful:
\begin{eqnarray}\label{polygonal}
S_N(t)&=& \frac{B_{\tau_1 +...+\tau_{\lfloor Nt \rfloor }} +
  \bigl(Nt-\lfloor Nt \rfloor \bigr) (B_{\tau_{\lceil
      Nt\rceil }} -B_{\tau_{\lfloor Nt\rfloor }})}{\sqrt{N}},\\ \label{scaled}
\BT_N(t)&=& \frac{B_{Nt}}{\sqrt{N}}.
\end{eqnarray}
Let $C([0,1],\mathbb{R}^k)$ be the space of continuous $k$-vector
valued functions on the unit interval equipped with the sup norm,
$\|\cdot \|_k$: for $f\in C([0,1],\mathbb{R}^k)$, $\|f\|_k
=\sum_{j=1}^k \sup_{t\in [0,1]} |f_j(t)|$.
Theorem~\ref{maintheorem} is proven in two steps. The first step
involves establishing probabilistic estimates on differences of
the above processes. The second step involves a Lipshitz bound on
the relevant Brownian concatenation operations and the application
of the quantitative estimates established in the first step.

Two estimates will be used in the proof of
Theorem~\ref{maintheorem}. Let $Y_i=\tau_i-1$ and $Z_i^N=\frac{1}{N}
\sum_{n=1}^i Y_n$. Note that $Z_i^N$ and $(Z_i^N)^2$ are submartingales. By
the Doob martingale inequalities \cite{RevuzYor},
\begin{eqnarray}
\PP \Bigr(\sup_{0\leq i\leq N} |Z_i^N|>\beta_N\Bigr) &\leq& \frac{\EE
  (Z_N^N)^2}{\beta_N^2} \nonumber \\
&=& \frac{\EE (\sum_{j=1}^N Y_i)^2}{N^2 \beta_N^2} = \frac{\EE \tau_1^2
  -(\EE \tau_1)^2}{N \beta_N^2} \nonumber \\
&\leq& \frac{Var(\tau_1)}{N \beta_N^2}.
\end{eqnarray}
For $\beta_N=\frac{c}{N^\lambda}$ where $c>0$ and $\lambda <\frac{1}{2}$, this inequality implies
\begin{equation}\label{estimate2}
\PP \Bigl( \sup_{1\leq i\leq N} \frac{|\sum_{j=1}^i\tau_j -i|}{N} >
\frac{c}{N^\lambda} \Bigr)\leq \frac{Var(\tau_1)}{N^{1-2\lambda}c^2}.
\end{equation}
This is the first estimate.

The second estimate concerns the modulus of continuity for Brownian
motion. If $\rho<\frac{\lambda}{2}$, then for some appropriately
chosen constants $A,\nu>0$,
\begin{equation}\label{estimate3}
\PP \Bigl( \sup_{t\in[0,1], |t-s|< cN^{-\lambda}} |\BT_N(t)-\BT_N(s)|
>\frac{c}{N^\rho} \Bigr) \leq Ae^{-\nu c}
\end{equation}
uniformly in $N$, $c>0$. This estimate is a consequence of Brownian
scaling and standard estimates for the maximum of Brownian motion.
The first part of the proof of Theorem~\ref{maintheorem} is complete.


Consider the functions $\Go:C([0,1],\mathbb{R}^k)\rightarrow
\mathbb{R}$ and $\Gu:C([0,1],\mathbb{R}^k)\rightarrow
\mathbb{R}$ defined by
\begin{eqnarray}
\Gu(f)=\inf_{0=t_0\leq t_1 \leq ... \leq t_k=1} \sum_{j=1}^k \bigl(f_j(t_j)
-f_j(t_{j-1})\bigr),\\
\Go(f)=\sup_{0=t_0\leq t_1 \leq ... \leq t_k=1} \sum_{j=1}^k \bigl(f_j(t_j)
-f_j(t_{j-1})\bigr).
\end{eqnarray}
\begin{lemma}\label{Lipshitz}
Both $\Go$ and $\Gu$ are Lipshitz with Lipshitz constant 2:
\begin{equation}
|\Go(f)-\Go(g)|\leq 2\|f-g\|_k \hspace{.1in} and \hspace{.1in} |\Gu(f)-\Gu(g)|\leq 2\|f-g\|_k,
\end{equation}
for all $f,g\in C([0,1],\mathbb{R}^k)$.
\end{lemma}
\noindent \textbf{proof}: Let $f, g\in C([0,1],\mathbb{R}^k)$ satisfy
$\|f-g\|_k=\epsilon$, $\|f_i-g_i\|_1=\epsilon_i$, and$\sum_{i=1}^k
\epsilon_i = \epsilon$ where $f=(f_1,...,f_k),
g=(g_1,...g_k)$. $|\Go(f)-\Go(g)| \leq \sum_{i=1}^k
|\Go(f_1,...,f_{i-1}, g_i,...,g_k)- \Go(f_1,...,f_i,g_{i+1},...g_k)|$
by the triangle inequality. $|\Go(f_1,...,f_{i-1}, g_i,...,g_k)-
\Go(f_1,...,f_i,g_{i+1},...g_k)|\leq 2\epsilon_i$. Thus,
$|\Go(f)-\Go(g)|\leq 2\|f-g\|_k$. Since the same argument holds for
$|\Gu(f)-\Gu(g)|$, the proof of the lemma is complete. $\Diamond$

By Skorohod embedding, for each $N>0$ construct a sequence of pairs of
random walks embedded in standard Brownian motions,
$\{(S_N^j(t),\BT_N^j(t))\}_{j=1}^\infty$, as in~\eqref{polygonal} and~\eqref{scaled}. Lemma~\ref{Lipshitz} and the estimates
on the differences of the Skorohod coupling imply that if $N,k\rightarrow
\infty$ and $k=o(N^\alpha )$, $\alpha< \frac{3}{14}$, then
\begin{eqnarray}\label{convergence}
k^\frac{1}{6}\bigl(\Go(S_N^1,...,S_N^k)
-\Go(\BT_N^1,...,\BT_N^k)\bigr)\rightarrow 0, \\
\label{convergence2} k^\frac{1}{6}\bigl(\Gu(S_N^1,...,S_N^k)
-\Gu(\BT_N^1,...,\BT_N^k)\bigr)\rightarrow 0,
\end{eqnarray}
in distribution. More explicitly, if $k=\lfloor N^\alpha \rfloor$,
then
\begin{eqnarray}
\PP \bigl( k^\frac{1}{6} \bigl(\Go(S_N^1,...,S_N^k)
-\Go(\BT_N^1,...,\BT_N^k)\bigr)>\epsilon\bigr) &\leq& \PP \Bigl(
 \sum_{i=1}^{\lfloor N^\alpha \rfloor}
\|S_N^i-\BT_N^i\|_1 >\frac{\epsilon}{2 N^\frac{\alpha}{6}}\Bigr)\nonumber \\
&\leq& \frac{N^\frac{7\alpha}{6}}{\epsilon} \EE \|S_N^1 - \BT_N^1\|_1.
\label{maininequality}
\end{eqnarray}
Note that
\begin{eqnarray}
\EE \|S_N^1-\BT_N^1\|_1 &=& \int_0^\infty \PP(\|S_N^1-\BT_N^1\|_1
>s)ds = \frac{1}{N^\rho} \int_0^\infty \PP \Bigl(\|S_N^1-\BT_N^1\|_1 >
\frac{u}{N^\rho}\Bigr)du \nonumber \\
&=& \frac{1}{N^\rho} \int_0^\infty \PP \Bigl(\|S_N^1-\BT_N^1\|_1 >
\frac{u}{N^\rho} ; \sup_{1\leq i \leq N} |\sum_{j=1}^i \tau_j -i|>
uN^{1-\lambda}\Bigr) du \nonumber \\
&& + \frac{1}{N^\rho} \int_0^\infty \PP \Bigl(\|S_N^1-\BT_N^1\|_1 >
\frac{u}{N^\rho} ; \sup_{1\leq i \leq N} |\sum_{j=1}^i \tau_j -i|<
uN^{1-\lambda}\Bigr) du.\label{final}
\end{eqnarray}
If $\lambda<\frac{1}{2}$, then~\eqref{estimate2}
implies that the integrand of the first integral on the right hand side
of~\eqref{final} has an integrable tail. If the
$\{X_i^j\}_{i,j=1}^\infty$ are Gaussian, this term is identically
$0$ since $\tau_i=1$ $\PP$-almost surely; in this case, there is no
constraint on $\lambda$.
The integrand of the second integral of~\eqref{final} is dominated by
\begin{equation}
\PP \Bigl( \sup_{t\in[0,1], |t-s|< uN^{-\lambda}} |\BT_N(t)-\BT_N(s)|
>\frac{u}{N^\rho} \Bigr).
\end{equation}
\eqref{estimate3} implies that if
$\rho<\frac{\lambda}{2}$, the second integrand
has exponentially decaying tails uniformly in $N$. A necessary condition for the
right hand side of~\eqref{maininequality} to vanish as $N\rightarrow
\infty$ is $\alpha< \frac{6}{7} \rho$. If $\{X_i^j\}_{i,j=1}^\infty$
are Gaussian, then $\rho< \frac{1}{2}$ is sufficient to make the
second integral of~\eqref{final} finite; thus, if
$\alpha<\frac{3}{7}$, the right hand side of~\eqref{maininequality} vanishes as
$N\rightarrow \infty$. If $\{X_i^j\}_{i,j=1}^\infty$ are not Gaussian,
then $\alpha< \frac{6}{7} \rho$ and $2\rho<\lambda<\frac{1}{2}$ imply
that $\alpha<\frac{3}{14}$ is sufficient to make~\eqref{final} vanish as $N\rightarrow \infty$.

As mentioned in the introduction, the theorems of
Baryshnikov \cite{Baryshnikov} and Gravener, Tracy, and Widom
\cite{GTW} state that $\Go(\BT_N^1,...,\BT_N^k)$ is distributed as
the largest eigenvalue of the GUE of $k\times k$ matrices.
By~\eqref{eq:GUElargest} and~\eqref{convergence},

\begin{equation}
k^\frac{1}{6}(\Go(S_N^1,...,S_N^k) - 2\sqrt{k})\Longrightarrow F_{GUE},
\end{equation}
in distribution. By~\eqref{convergence2} and the symmetry with
respect to multiplication by $-1$ of Brownian
motion and the Dyson process
\begin{equation}
k^\frac{1}{6}(-\Gu(S_N^1,...,S_N^k) - 2\sqrt{k})\Longrightarrow
F_{GUE}.
\end{equation}
Since $\Go(S_N^1,...,S_N^k)$ and $\Gu(S_N^1,...,S_N^k)$
are distributed as $L(N,k)$ and $R(N,k)$, respectively, the proof of
Theorem~\ref{maintheorem} is complete.
\end{proof}

\medskip

It is clear from the proof of Theorem~\ref{maintheorem} that when
$k>>o(N^{\frac37})$  the `up movements' of the paths contribute to
$L^l(N,k)$. Hence, one needs to take those terms into account in order
to prove a central limit theorem for general $N,k$ scaling. The
arguments in this paper do not seem to extend to the general case.

The proof of Corollary~\ref{LPP} remains.
\begin{proof}[proof of Corollary~\ref{LPP}]
With no loss of generality, assume that $\EE X_i^j=0$ and $\EE |X_i^j|^2=1$. The definitions of $L^l(N,k)$ and $L(N,k)$ imply the following
deterministic estimate:
\begin{equation}
k^\frac{1}{6}\inf_{1\leq i_1 \leq \cdots \leq i_{k-1}\leq N} \Big \{ \frac{1}{\sqrt{N}}
\sum_{j=1}^{k-1} X_{i_j}^{j+1} \Big \}\leq
k^\frac{1}{6}\frac{L^l(N,k)-L(N,k)}{\sqrt{N}} \leq
k^\frac{1}{6}\sup_{1\leq i_1 \leq \cdots \leq i_{k-1}\leq N} \Big \{\frac{1}{\sqrt{N}}
\sum_{j=1}^{k-1} X_{i_j}^{j+1}\Big \}.
\end{equation}
Note that
\begin{eqnarray} \label{CorollaryEstimate}
\PP \Bigl( k^\frac{1}{6}\sup_{1 \leq i_1 \leq \cdots \leq i_{k-1}\leq N} \frac{1}{\sqrt{N}}
\sum_{j=1}^{k-1} X_{i_j}^{j+1} >\epsilon \Bigr) &\leq&
\frac{k^\frac{1}{6}}{\epsilon \sqrt{N}} \EE \Bigl(\sup_{1 \leq i_1 \leq
  \cdots \leq i_{k-1}\leq N}\sum_{j=1}^{k-1} X_{i_j}^{j+1}\Bigr)  \nonumber \\
&\leq&
\frac{k^\frac{7}{6}}{\epsilon \sqrt{N}} \EE \Bigl(\sup_{i\in\{1,...,N\}}
|X_i^1|\Bigr) \nonumber \\
&\leq&  \frac{k^\frac{7}{6}}{\epsilon \sqrt{N}} \int_0^\infty \Bigl \{
1- \PP \Bigl( |X_1^1| <s \Bigr)^N \Bigr \} ds \nonumber \\
&\leq& \frac{k^\frac{7}{6}}{\epsilon \sqrt{N}} \int_0^\infty \Bigl \{ 1-\Bigl [ 1- \PP \Bigl( |X_1^1| >s \Bigr) \Bigr ]^N \Bigr \}ds \nonumber \\
&=& \frac{k^\frac{7}{6}}{\epsilon N^{\frac{1}{4}}} \int_0^\infty \Bigl \{
1-\Bigl [ 1- \PP \Bigl( |X_1^1| >s N^{\frac14} \Bigr) \Bigr ]^N
\Bigr \}ds.
\end{eqnarray}
Now split the last integral into two pieces; one over the interval
$[0,1]$ and the other over $(1,\infty)$. Using the basic estimate
$0\le \PP \Bigl( |X_1^1| >s N^{\frac14} \Bigr) \le 1$ for the
first integral and the Chebyshev inequality for the second
integral,~\eqref{CorollaryEstimate} is less than or equal to
\begin{equation}\label{CorollaryEstimate2}
   \frac{k^\frac{7}{6}}{\epsilon N^{\frac14}} \int_0^1 ds
+ \frac{k^\frac{7}{6}}{N^{\frac14}} \int_1^\infty \Bigl \{ 1-\Bigl
[ 1- \frac{\EE |X_1^1|^4}{s^4 N} \Bigr ]^N \Bigr \}ds \le
\frac{k^\frac{7}{6}}{\epsilon N^{\frac14}} +
\frac{k^\frac{7}{6}}{\epsilon N^{\frac14}} \int_1^\infty \Bigl \{
1- e^{-2\frac{\EE |X_1^1|^4}{s^4}} \Bigr \}ds.
\end{equation}
If  $N,k\rightarrow \infty$ such that $k=o(N^\alpha)$, $\alpha<
\frac{3}{14}$, then the right hand side
of~\eqref{CorollaryEstimate2} tends to $0$. The $\inf$ can be
treated similarly. Theorem~\ref{maintheorem}
and~\eqref{CorollaryEstimate},~\eqref{CorollaryEstimate2} imply
Corollary~\ref{LPP} for the non-Gaussian case. If
$\{X_i^j\}_{i,j=1}^\infty$ are Gaussian, then the arguments from
the proof of Theorem~\ref{maintheorem} apply.
\end{proof}



\end{document}